\documentclass{article}
\usepackage[left=1.1in, right=1.1in, top=0.9in, bottom=1.1in]{geometry}
\usepackage{amsmath}
\usepackage{amsfonts}
\usepackage{amssymb}
\usepackage{amsxtra}
\usepackage{amsthm}
\usepackage{mathrsfs}
\usepackage{color}
\usepackage{enumerate}
\usepackage{graphicx}
\usepackage{float}

\graphicspath{ {graph/} }

\newtheorem*{theorem*}{Theorem}
\newtheorem*{corollary*}{Corollary}
\newtheorem*{definition*}{Definition}
\newtheorem*{remark*}{Remark}
\newtheorem*{proposition*}{Proposition}

\begin{document}
\title{Learning Algorithms for Coarsening Uncertainty Space and Applications to Multiscale Simulations}
\author{Zecheng Zhang\thanks{Department of Mathematics, Texas A\&M University,
College Station, TX, USA (\texttt{tom\_z\_z@tamu.edu})} ,
Eric Chung\thanks{Department of Mathematics, The Chinese University of Hong Kong, Shatin, New Territories, Hong Kong SAR, China (\texttt{tschung@math.cuhk.edu.hk}) },
Yalchin Efendiev\thanks{Department of Mathematics \& Institute for Scientific Computation (ISC), Texas A\&M University,
College Station, TX, USA (\texttt{efendiev@math.tamu.edu}) \& Multiscale Mathematics Laboratory, North-Eastern Federal University, Russia },
Wing Tat Leung\thanks{ICES, University of Texas, Austin, TX, USA (\texttt{wleung@ices.utexas.edu})}
}

\maketitle

\begin{abstract}
In this paper, we investigate and design multiscale simulations for 
stochastic
multiscale PDEs. As for the space, we consider a coarse grid and a known
multiscale method, the Generalized Multiscale Finite Element Method (GMsFEM). 
In order to obtain a small dimensional representation of the solution
in each coarse block, the uncertainty space needs to be partitioned 
(coarsened). This coarsenining collects realizations that provide similar
multiscale features as outlined in GMsFEM (or other method of choice).
This step is known to be computationally demanding as it requires
many local solves and clustering based on them (see \cite{paper}).
In this paper, we take a different approach and learn coarsening
the uncertainty space. Our methods use deep learning techniques in 
identifying clusters (coarsening) in the uncertainty space. 
We use convolutional neural networks combined with some techniques in
adversary neural networks. We define appropriate loss functions 
in the proposed neural networks, where the loss function is composed
of several parts that includes terms related to clusters and reconstruction
of basis functions. We present numerical results for channelized permeability
fields in the examples of flows in porous media.
\end{abstract}

\section{Introduction}


Many problems are multiscale with uncertainties. Examples
include problems in 
porous media, material sciences, biological sciences, and so 
on. For example, in porous media applications, engineers 
can obtain fine-scale
data about pore geometries or subsurface properties at very fine resolutions. 
These data are  
obtained in some spatial locations and then generalized to the entire
reservoir domain. As a result, one uses geostatistical or other statistical
tools to populate the media properties in space. The resulting porous media
properties are stochastic and one needs to deal with many porous media 
realizations, where each realization is multiscale and varies at
very fine scales.

Simulating each realization can be computationally
expensive because of the media`s multiscale nature. 
Our objective is to simulate  many of these realizations.
To address the issues associated with spatial and temporal
scales,  many multiscale methods have been developed
(\cite{gmsfem2, gmsfem3, gmsfem4, gmsfem5, gmsfem6, gmsfem7, gmsfem8,gmsfem9,ye1, ElasticGMsFEM, chung2017DGstokes}). 
These methods perform simulations on the coarse grid
by developing reduced-order models. However, developing reduced-order
models requires local computations, which can be expensive when
one deals with many realizations. For this reason, some type of
coarsening of the uncertainty space is needed (\cite{paper}).
In this paper, we consider some novel approaches for developing
coarsening of uncertainty space as discussed below.

To coarsen the uncertainty space, clustering algorithms are often used;
but a proper distance function should be designed 
in order to make the clusters have physical sense and achieve a reduction
in the uncertainty space.
The paper (\cite{paper}) proposed a method that uses the distance 
between local solutions. The motivation is that 
the local problems with random boundary conditions can 
represent the main models with all boundary conditions.
Due to a high dimension of the uncertainty space, the authors in 
(\cite{paper}) proposed to compute the local solutions of only several realizations and then use the 
Karhunen-Loeve expansion(\cite{kl1}) to approximate the solutions of all the other realizations. The distance function is then defined to be the distance between solutions and the standard K-means (\cite{kmeans}) algorithm is used to cluster the uncertainty space.

The issue with this method is computing the local solutions in the local neighborhoods. It is computationally expensive to compute the local solutions; although the KL expansion can save time to approximate the solutions of other realizations, one still needs to decide how many selected realizations we need to represent all the other solutions. In this paper, we propose the use of deep learning methodology and avoid explicit clustering as in earlier works.
We remark that the development of deep learning techniques for multiscale simulations are recently reported in \cite{wang2020deep,wang2020reduced,vasilyeva2019learning,cheung2019deep,wang2019prediction}.

In this work, to coarsen the uncertainty space, we propose a deep learning algorithm which will learn the clusters for each local neighborhood. Due the nature of the permeability fields, we can use the transfer learning which uses the parameters of one local neighborhood to initialize the learning of all the other neighborhoods. This saves significantly computational time.

The auto encoder structure \cite{goodfellow} has been widely used in improving the K-mean clustering algorithm (\cite{cluster1,cluster2, cluster3}).
The idea is to use the encoder to extract features and reduce the dimension; the encoding process can also be taken as a kernel method \cite{bishop} which maps the data to a space which is easier to be separated. The decoder is used to upsample the latent space (reduced dimension feature space) back to the input space. The clustering algorithm is then used to cluster the latent space, which will save time due to the low dimension of the latent space and also preserve the accuracy due to the features extracted by the encoder.

Traditionally, the learning process is only involved in reconstructing the input space. Such kind of methods ignore the features extracted by latent space; so, it is not clear if the latent space is good enough to represent the input space and is easily clustered by the K-means method. In \cite{cluster3}, the authors proposed a new loss which includes the reconstruction loss meanwhile the loss results from the clustering. The authors claimed that the new loss improves the clustering results.

We will apply the auto encoder structure and the multiple loss function; 
however, we will design the auto encoder as a generative network, i.e., the input and output space are different. 
More precisely, the input is the uncertain space (permeability fields) and the output will be the multiscale functions co-responding to the uncertain space. 
Intuitively, we want to use the multiscale basis to supervise the learning of the clusters so that the clusters will inherit the property of the solution.
The motivation is the multiscale basis can somehow represent the real solutions and permeability fields;
hence, the latent space is no longer good for clustering the input space but will be suitable for representing the multiscale basis function space.

To define the reconstructing loss, the common idea is the mean square error (MSE); but many works (\cite{condition, supergan, laplacian, perceptual}) have shown that the MSE tends to produce the average effect. In fact, in the area of image super-resolution
(\cite{condition,supergan, laplacian, perceptual, dong, kim1, residualchannel, residualdense, recursive, enhanced, memnet}) and other low level computer vision tasks, the generated images are usually over-smooth if trained using MSE. The theory is the MSE will capture the low frequency features like the background which is relatively steady; but for images with high contrast, the MSE will usually try to blur the images and the resulting images will lose the colorfulness and become less vivid (\cite{condition}). Our problem has multiscale nature and we want to capture the dominant modes and multiscale features, hence a single MSE is clearly not enough.

Following the idea from (\cite{perceptual, supergan}), we consider adding an adversary net (\cite{gan}). The motivation is the fact that different layers of fully convolutional network extract different features (\cite{senet, perceptual, selfattentiongan}). 
Deep fully convolutional neural networks (FCN) (\cite{long,inception, rethinking, segnet, deconvolution, dense}) have demonstrated its power in almost all computer vision tasks. Convolution operation is a local operation and the network with full convolutions are independent with the input size. People now are clear about the functioning of the different layers of the FCN.
In computer vision task, the lower layers (layers near input) tend to generate sharing features of all objects like edges and curves while the higher layers (near output) are more object oriented. If we train the network using the loss from the lower layers, the texture and details are persevered, while the higher layers will keep the general spatial structure. 

This motivates us using the losses from different layers of the fully convolutional layers. Multiple layers will give us a multilevel capture of the basis features and hence measure the basis in a more complete way. To implement the idea, we will pretrain an adversary net; and input the multiscale basis of the generative net and the real basis. The losses then come from some selected layers of the adversary net. 
Although it is still not clear the speciality of each layer, if we consider the multiscale physical problem, the experiments show that the accuracy is improved and, amazingly, the training becomes easier when compared to the MSE of the basis directly.

The uncertain space coarsening (cluster) is performed using the deep learning idea described above. But due to the space dimension, we will perform the clustering algorithm locally in space; that is, we first need a spatial coarsening. 
Due to the multiscale natural of the problem, 
this motivates us using the generalized multiscale finite element methods (GMsFEM) which derives the multiscale basis of a coarse neighborhood by solving the local problem. 
GMeFEM was first proposed in (\cite{gmsfem1}) and further studied in (\cite{gmsfem2, gmsfem3, gmsfem4, gmsfem5, gmsfem6, gmsfem7, gmsfem8, gmsfem9}). This method is
a generalization of the multiscale finite element method (\cite{msfem1, msfem2}). The work starts from constructing the snapshot space for each local neighborhood. The snapshot space is constructed by solving local problems and several methods including harmonic extension, random boundary condition \cite{snapshot} have been proposed. Once we have the snapshot space, the offline space which will be used as computing the solution are constructed by using spectral decomposition. 

The rest of the work is organized as follow: in Section 2, we 
consider the problem setup
and introduce both uncertain space and spatial coarsening.
In Section 3, we introduce the structure of the network and the training algorithm. In Section 4, we will present the numerical results.
The paper ends with conclusions. 

\section{Problem Settings}

In this section, we will present some basic ideas involving the use of the generalized multiscale finite element method (GMsFEM)
for parameter-dependent problems. 
Let $D$ be a bounded domain in $\mathbb{R}^2$ and $\Omega$ be the parameter space in $\mathbb{R}^N$. We
 consider the following parameter-dependent elliptic problem:
\begin{align}
-\nabla\cdot(\kappa(x,s)\nabla u(x,s)) = f(x,s), (x,s)\in D\times\Omega \label{eqn1},\\
u(x,s) = 0, (x,s)\in\partial D\times\Omega,
\label{eqn2}
\end{align}
where $\kappa(x,s)$ is a heterogeneous coefficient depending on both the spatial variable $x$ and the parameter $s$, and $f\in L^2(D)$ is a given source. 
We remark that the differential operators in (\ref{eqn1}) are defined with respect to the spatial variable $x$.
This is the case for the rest of the paper. 

\subsection{The Coarsening of the Parameter Space. The Main Idea.}

The parameter space $\Omega$ is assumed to be of very high dimension (i.e. large $N$) and consists of 
very large number of realizations.
For a given realization, the idea is to find its representation in the 
coarse space and use the coarse space to perform the computation. 
We will use the deep cluster learning algorithm to perform the coarsening. 
Due to the heterogeneous properties of the proposed problem, 
fine mesh is used; this will bring difficulties in coarsening 
the parameter space and in computation of the solution. 
We hence perform the parameter coarsening locally in the space $D$, 
i.e., we also coarsen the spatial domain. 
To coarsen the spatial domain, we use coarse grids and 
consider the GMsFEM.

\begin{figure}[H]
\centering
\includegraphics[scale = 0.4]{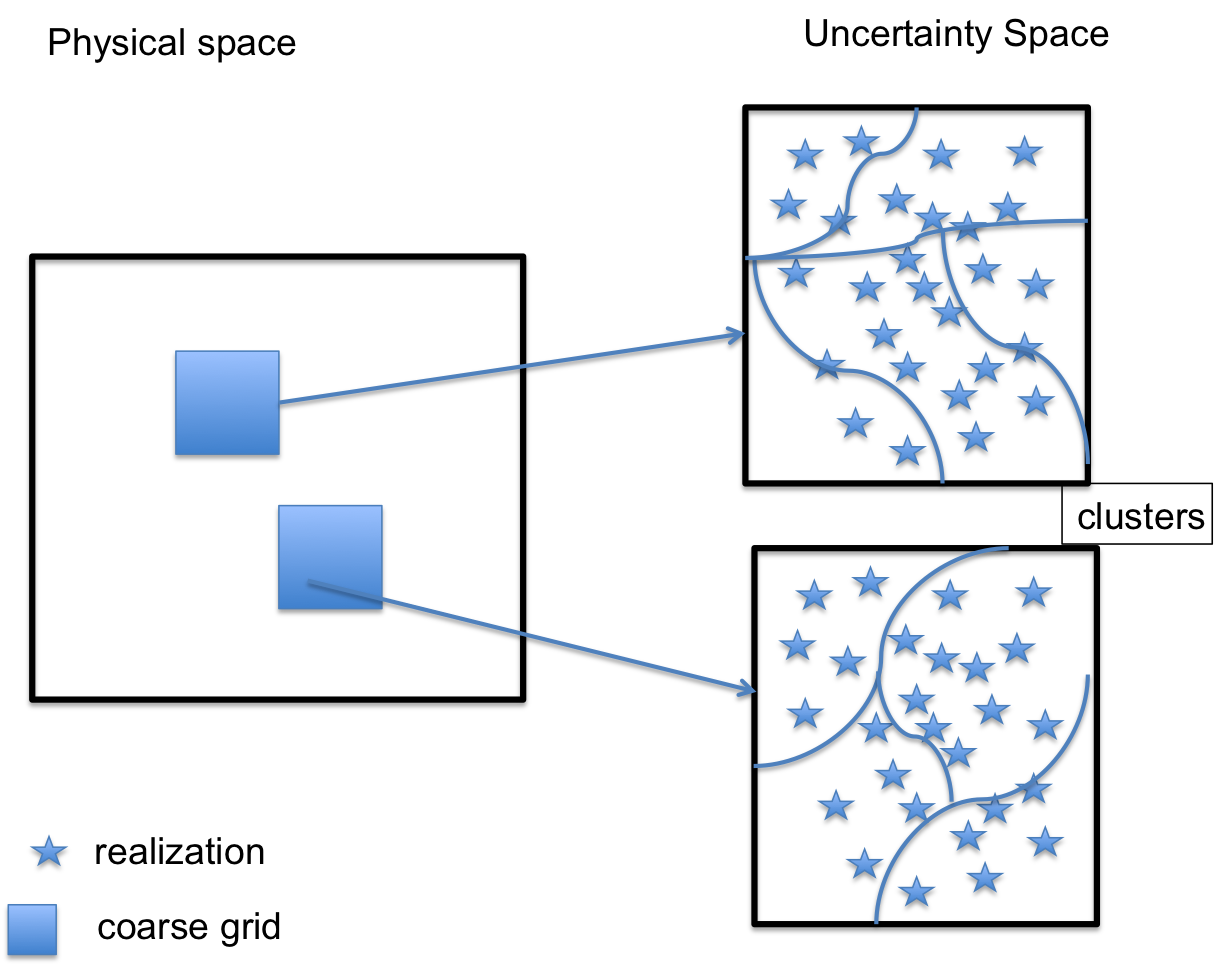}
\caption{Illustration of coarsening of space and uncertainties. Different clusters for different coarse blocks. On the left plot, two coarse blocks are shown. On the right plot, clusters are illustrated. }
\label{ill1}
\end{figure}

In Figure \ref{ill1}, we present an illustration of the proposed 
coarsening technique. 
On the left figure, the coarse grid blocks in the space are shown.
Each coarse grid has a different clusters in the uncertainty space $\Omega$,
which corresponds to the coarsening of the uncertainty space. The main
objective in multiscale methods is efficiently finding the clustering
of the uncertainty space, which is our main goal. 

\subsection{Space Coarsening --- Generalized Multiscale Finite Element Method}

It is computationally expensive to capture heterogeneous properties using very fine mesh. For this reason, we use GMsFEM to coarsen the spatial representation of the solution. The coarsening of the parameter space will be performed in each local spatial neighborhood. We will achieve this goal by the GMsFEM, which will briefly be discussed. 
Consider the second order elliptic equations $Lu = f$ in $D$ with proper boundary conditions; denote the the elliptic operator as:
\begin{align}
L(u) = -\frac{\partial}{\partial x_i}(k_{ij}(x)\frac{\partial}{\partial x_j}u).    
\end{align}

\begin{figure}[H]
\centering
\includegraphics[scale = 0.2]{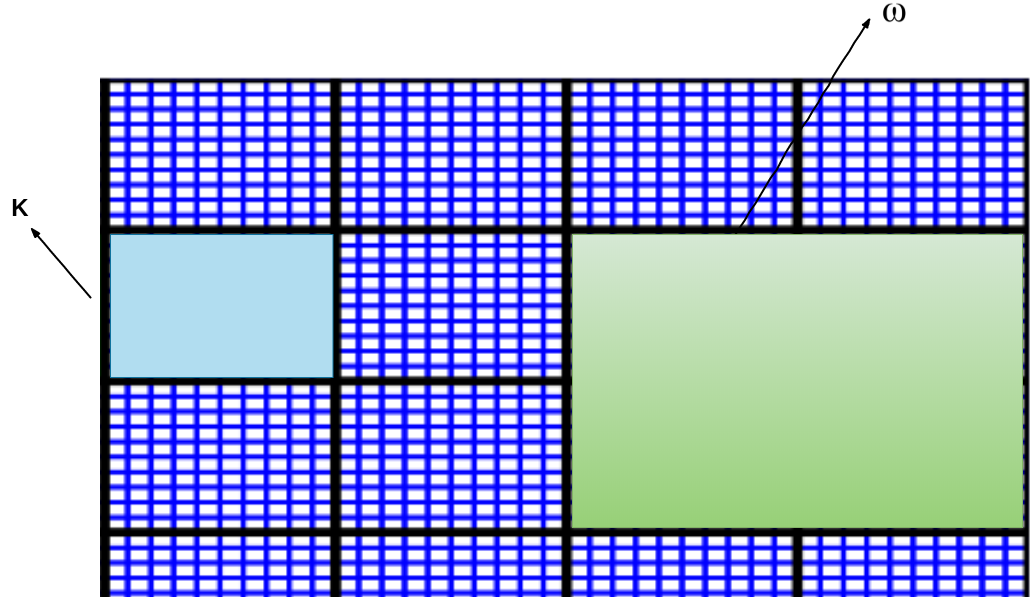}
\caption{Domain Partition $\mathcal{T}^H$.}
\label{grid}
\end{figure}

Let the spatial domain $D$ be partitioned by a coarse grid $\mathcal{T}^H$; this does not resolve the multiscale features. Let us denote $K$ as one cell in $\mathcal{T}^H$ and refine $K$ to obtain the fine grid partition $\mathcal{T}^h$ (blue box in Figure~\ref{grid}). We assume the fine grid is a conforming refinement of the coarse grid. 
See Figure~\ref{grid} for details. 

For the $i$-th coarse grid node, let $\omega_i$ be the set of all coarse elements having the vertex $i$ (green region in Figure~\ref{grid}).
We will solve local problem in each coarse neighborhood to obtain set of multiscale basis functions $\{ \phi_i^{\omega_i}\}$
and seek solution in the form:
\begin{align}
    u = \sum_i\sum_j c_{ij}\phi_j^{\omega_i},
\end{align}
where $\phi_j^{\omega_i}$ is the offline basis function in the $i$-th coarse neighborhood $\omega_i$ and $j$ denotes the $j$-th basis function. 
Before we construct the offline basis, we first need to derive the snapshot basis.

\subsubsection{Snapshot Space }
There are several methods to construct the snapshot space; we will use the harmonic extension of the fine grid functions defined on the boundary of $\omega_i$.  Let us denote $\delta_l^h(x)$ as fine grid delta function, which is defined as $\delta_l^h(x_k) = \delta_{l,k}$ for $x_k\in J_h(\omega_i)$ where $J_h(\omega_i)$ denotes the boundary nodes of $\omega_i$. The snapshot function $\psi_l^{\omega_i}$ is then calculated by solving local problem in $\omega_i$:
\begin{align}
    L(\psi_l^{\omega_i}) = 0
\end{align}
subject to the boundary condition $\psi_l^{\omega_i} = \delta_l^h(x)$. The snapshot space $V_{snap}^{\omega_i}$ is then constructed as the span of all snapshot functions.

\subsubsection{Offline Spaces}
The offline space $V_{off}^{\omega_i}$ is derived from the snapshot space and is used for computing the solution of the problem. We need to solve for a spetral problem and this can be summarized as finding $\lambda$ and $v\in V_{snap}^{\omega_i}$ such that:
\begin{align}
    a_{\omega_i}(v,w) = \lambda s_{\omega_i}(v,w), \forall w\in V_{snap}^{\omega_i},
    \label{eig}
\end{align}
where $a_{\omega_i}$ is symmetric non-negative definite bilinear form and $s_{\omega_i}$ is symmetric positive definite bilinear form. By convergence analysis, they are given by
\begin{align}
    &a_{\omega_i}(v,w) = \int_{\omega_i}\kappa \nabla v\cdot \nabla w,\\
    &s_{\omega_i}(v,w) = \int_{\omega_i} \tilde{\kappa} v\cdot w.
\end{align}
In the above definition of $s_{\omega_i}$, the function $\tilde{\kappa} = \kappa \sum |\nabla \chi_j|^2$
where $\{ \chi_j\}$ is a set of partition of unity functions corresponding to the coarse grid partition of the domain $D$ 
and the summation is taken over all the functions in this set. 
The offline space is then constructed by choosing the smallest $L_i$ eigenvalues and we can form the space by the linear combination of snapshot basis using corresponding eigenvectors:
\begin{align}
    \phi_k^{\omega_i} = \sum_{j = 1}^{L_i} \Psi_{kj}^{\omega_i}\psi_j^{\omega_i},
\end{align}
where $\Psi_{kj}^{\omega_i}$ is the $jth$ element of $kth$ eigenvector and $L_i$ is the number of snapshot basis. $V_{off}$ is then defined as the collection of all local offline basis functions. Finally we are trying to find $u_{off}\in V_{off}$ such that
\begin{align}
    a(u_{off},v) = \int_{D}fv, \forall v\in V_{off}
    \label{linear}
\end{align}
where $a(u,v) = \int_D \kappa \nabla u \cdot \nabla v$.
For more details, we refer the readers to the references \cite{gmsfem7,gmsfem8,gmsfem9}.

\subsection{The Idea of the Proposed Method}

We present the general methodology in this section. The target is to save the time in computing the GMsFEM basis $\phi_k^{\omega_i}$ for all $\omega_i$ and for all uncertain space parameters.
We propose the clustering algorithm to coarsen the uncertain space in each local neighborhood. The key to the success of the clustering is that: the cluster should inherit the property of the solution, that is, the local heterogeneous fields $\kappa(x,s)$ clustered into the same group should have similar solution properties. When the cluster is learnt by the some learning algorithm, the only computation involved is to fit the local neighborhood of the given testing heterogeneous field into some cluster. This is a feed forward  process including several convolution operations and matrix multiplications and compared to the direct computing, we save a lot of time in computing the spectral problem (\ref{eig}) and the inverse of a matrix (\ref{linear}). The detailed process is illustrated in the following chart (Figure~\ref{workflow}):
\begin{figure}[H]
\centering
\includegraphics[scale = 0.5]{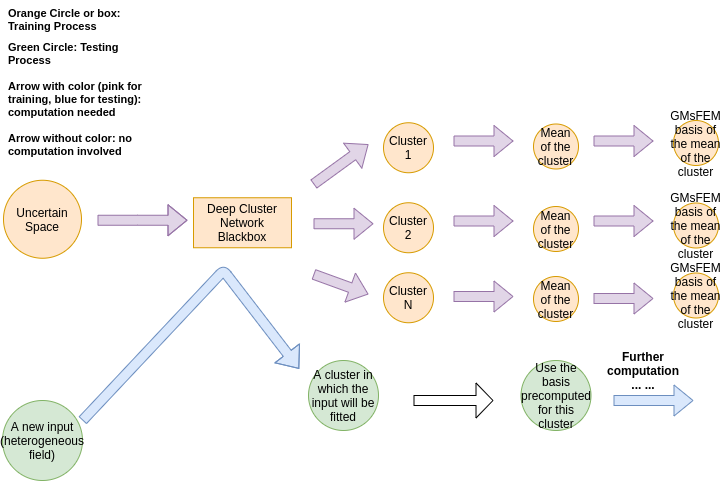}
\caption{Work Flow of the Proposed Method}
\label{workflow}
\end{figure}
\begin{enumerate}
\item (Training) For a given input local neighborhood $\omega_j$, we train the cluster (which will be detailed in next section) of the parameter space $\Omega$
and get the clusters $S_1^j,..., S_n^j$, where $n$ is the number of clusters and is uniform for all $j$. Please note that we may have different cluster assignments in different local neighborhoods.
\item (Training) For each local neighborhood $\omega_j$ and cluster $S_i^j$, define the average $\Bar{\kappa}_{ij}$ and compute generalized multiscale basis for $\Bar{\kappa}_{ij}$.

\item (Testing) Given a new $\kappa(x,s)$ and for each local neighborhood $\omega_j$, fit $\kappa(x,s)$ into a $S_i^j$ by the trained network (step 1) and use the pre-computed GMsFEM basis (step 2) to find the solution.
\end{enumerate}
It should be noted that we perform clustering using the heterogeneous fields; however, the cluster should inherit the property of the solution corresponding to the heterogeneous fields. This makes the clustering challenging. The performance of the standard K-Means algorithm relies on the initialization and the distance metric. We may initialize the algorithm based on the clustering of the heterogeneous fields but we need to design a good metric. In the next section, we are going to introduce a learning algorithm which uses an auto-encoder structure and multiple losses to achieve the required clustering task.

\section{Deep Learning}
The network is consisted of two sub networks. The first one is targeted to performing the clustering and second one, which is the adversary net, will serve as the reconstruction of loss function.
\subsection{Clustering Net}
\begin{figure}[H]
\centering
\includegraphics[scale = 0.4]{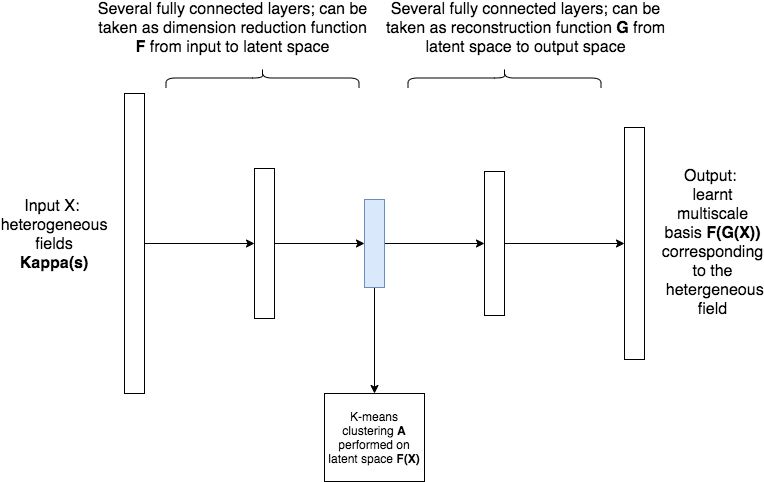}
\caption{cluster Network}
\label{clusternet}
\end{figure}

The cluster net is aimed for clustering the heterogeneous fields $\kappa(x,s)$; but the resulting clusters should inherit the properties of the solution corresponding to the $\kappa(x,s)$, i.e., the heterogeneous fields grouped in the same cluster should have similar corresponding solution properties. This similarity will be measured by the adversary net which will be introduced in section (\ref{advnet}).
We hence design the network demonstrated in Figure~\ref{clusternet}.
The input $X\in\mathbb{R}^{m,d}$ ,where $m$ is the number of samples and $d$ is the dimension of one local heterogeneous field, of the network is the local heterogeneous fields which are parametrized by the random variable $s\in\Omega$. The output of the network is the multiscale basis (first GMsFEM basis) which somehow represents the solution corresponding to the coefficient $\kappa(x,s)$. This is a generative network which has an auto encoder structure. The dimension reduction function $F(X)$ can be interpreted as some kind of kernel method which maps the input data to a new space which is easier to be separated. This can also be interpreted as the learning of a good metric function for the later performed K-mean clustering.
We will perform K-means clustering algorithm in latent space $F(X)$. 
$G(\cdot)$ will then transfer the latent space data to the space of multiscale basis function. This process can be taken as a generative process and we reconstruct the basis from the extracted features.
The detailed algorithm is as follow (see Figure~\ref{deepalgo} for an illustration):
\begin{figure}[H]
\centering
\includegraphics[scale = 0.4]{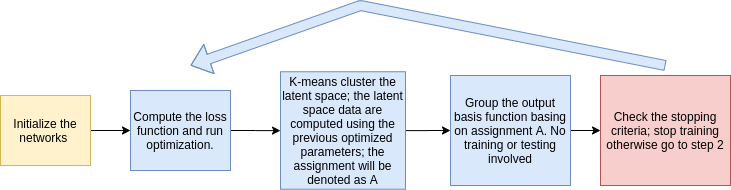}
\caption{Deep Learning Algorithm }
\label{deepalgo}
\end{figure}

\noindent
{\bf Steps illustrated in Figure~\ref{deepalgo}}:
\begin{enumerate}
\item Initialize the networks and clustering the output basis function.
\item
Compute the loss function $L$ (defined later) and run optimization.
\item 
Cluster the latent space by K-Means algorithm (reduced dimension space, which is a middle layer of the cluster network); the latent space data are computed using the previous optimized parameters; the assignment will be denoted as $A$.
\item 
Basis functions whose corresponding inputs are in the same cluster (basing on assignment A) will be grouped together. No training or fitting-in involved in this step. 
\item 
Repeat step $2$ to step $4$ until the stopping criteria is met.
\end{enumerate}

\subsection{Loss Functions}
Loss function is the key to the deep learning. Our loss function is consisted of cluster loss and the reconstruction loss.
\begin{enumerate}
\item
Clustering loss $C(\theta_F,\theta_G)$: this is the mean standard deviation of all clusters of the learnt basis and $\theta$ is the parameters we need to optimize.
It should be noted that the loss here is computed using the learnt basis instead of the input of the network.
This loss controls the clustering process, i.e., the smaller the loss, the better the clustering in the sense of clustering the multiscale basis. 
Let us denote $\kappa_{ij}$ as $jth$ realization in $ith$ cluster;   $G(F(\kappa_{ij}))\in\mathbb{R}^d$ will then be $jth$ learnt basis in cluster $i$ and let $\theta_G$ and $\theta_F$ be the parameters associated with $G$ and $F$, the loss is then defined as follow,
\begin{align}
C(\theta_F,\theta_G) = \frac{1}{|A|}\sum_{i}^{|A|}\sum_j^{A_i} \frac{1}{A_i}\|G(F(\kappa_{ij}))- \Bar{\phi}_i\|_2^2,
\end{align}
where $|A|$ is the number of clusters which is a hyper parameter and $A_i$ denotes the number of elements in cluster $i$; $\Bar{\phi}_i\in\mathbb{R}^d$ is the mean of cluster $i$. 
This loss clearly serves the purpose of clustering the solution instead of the input heterogeneous fields; however, in order to guarantee the learnt basis are closed to the pre-computed multiscale basis, we need to define the reconstruction loss which measures the difference between the learnt basis and the pre-computed basis. 

\item Reconstruction loss $R(\theta_F,\theta_G)$: this is the mean square error of multiscale basis $Y\in\mathbb{R}^{m,d}$, where $m$ is the number of samples.
This loss controls the construction process, i.e., if the loss is small, the learnt basis are close to the real multiscale basis. This loss will supervise the learning of the cluster. It is defined as follow:
\begin{align}
    R(\theta_F,\theta_G) = \frac{1}{m}\sum_{i}^m \|G(F(\kappa_{i}))- \phi_{i}\|_2^2, \label{recon_err}
\end{align}
where $G(F(\kappa_{i}))\in\mathbb{R}^d$ and $\phi_i\in \mathbb{R}^d$ are learnt and pre-computed multiscale basis of $ith$ sample $\kappa_i$ separately. 
\end{enumerate}
The entire loss function is then defined as $L(\theta_F,\theta_G) = \lambda_1 C +\lambda_2 R$, where $\lambda_1, \lambda_2$ are predefined weights. We are going to solve the following optimization problem:
\begin{align}
    \min_{\theta_G,\theta_F} L(\theta_F,\theta_G)
\end{align}
for the required training process. 

\subsection{Adversary Network Severing as an Additional Loss \label{advnet}}
We have introduced the reconstruction loss which measures the similarity between the learnt basis and the pre-computed basis in the previous section. It is the Mean square error (MSE) of the learnt and pre-computed basis. 
MSE is a smooth loss and easy to train but there is a well known fact about MSE that this loss will blur the image.
In the area of image super-resolution and other low level computer vision tasks, the loss is not friendly to inputs with high contrast and the resulting generated images are usually over smooth. Our problem has multiscale nature and is similar with the low level computer vision task, i.e., this is a generative task; hence blurring and over smooth should happen if the model is trained by MSE. To define a great reconstruction loss is important.

Motivated by some works about the successful application of deep fully convolutional network (FCN) in computer vision, we design a perceptual loss to measure the error. It is now clear that the lower layers in the FCN usually will extract some general features shared by all objects like the horizontal (vertical) curves, while the higher layers are usually more objects oriented. This gives people the insight to train the network using different layers. Johnson then proposed the perceptual loss \cite{perceptual} which is the combination of the MSE of selected layers of the VGG model \cite{simonyan2014very}. The authors claim in their paper that the early layers tends to produce images that are visually indistinguishable from the input; however if reconstruct from higher layers, image content and overall spatial structure are preserved but color, texture, and exact shape are not.

\begin{figure}[H]
\centering
\includegraphics[scale = 0.4]{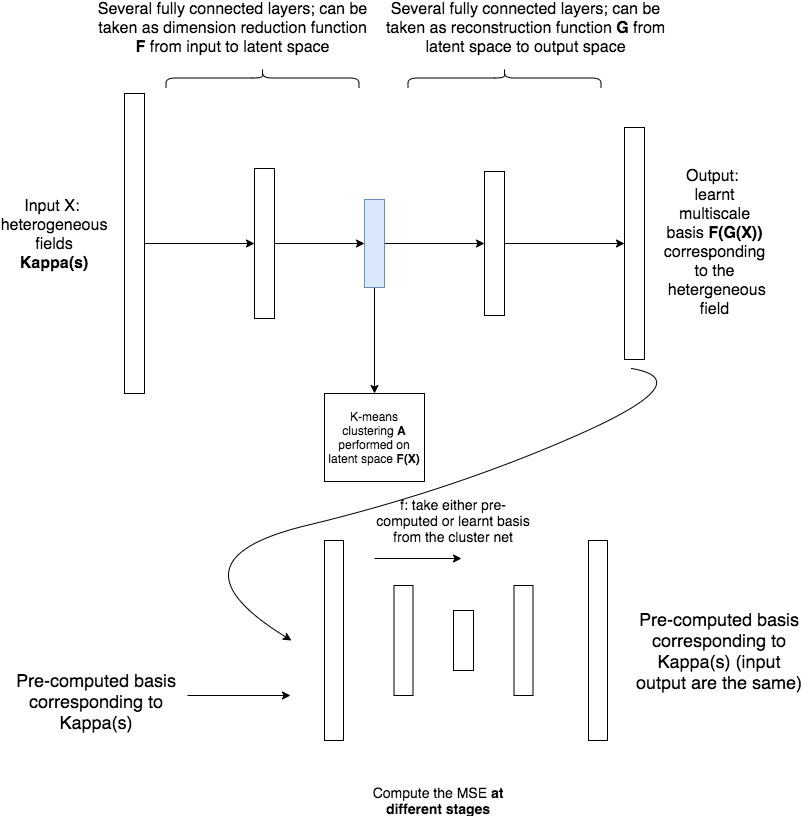}
\caption{The Complete Network}
\label{entirenet}
\end{figure}

We will adopt the perceptual loss idea and design an adversary network to compute an additional reconstruction loss. The network structure can be seen in Figure~\ref{entirenet}.
The adversary net is fully convolutional with input and output both pre-computed multiscale basis. The network has an auto encoder structure and is pre-trained; i.e., we are going to solve the following minimization problem:
\begin{align}
\min_{\theta_A} \frac{1}{m}\sum_{i}\|f(\phi_i)-\phi_i)\|_2^2,
\end{align}
where $\phi_i$ is the multiscale basis and $f$ is the adversary net associated with trainable parameter $\theta_A$.
Denote $f_j(\cdot)$ as the output of layer $j$ of the adversary network. The additional reconstruction loss is then redefined as:
\begin{align}
    A(\theta_F,\theta_G) = \frac{1}{m}\sum_{i = 1}^m\sum_{j\in I} \|f_j(G(F(\kappa_{i})))-f_j(\phi_i)\|^2_2, \label{adv_err}
\end{align}
where $I$ is the index set which contains some layers of the adversary net. The complete optimization problem can be now formulated as:
\begin{align}
   \min_{\theta_G,\theta_F}\lambda_1 C +\lambda_2 R+\lambda_3 A. \label{loss_all}
\end{align}

\section{Numerical Experiments}
In this section, we will demonstrate a series of experiments. We are going to apply our method on problems with high contrast including moving background and moving channels. Let us first introduce the high contrast field.

\subsection{High Contrast Heterogeneous Fields with Moving Channels}
\begin{figure}[H]
\centering
\includegraphics[scale = 0.5]{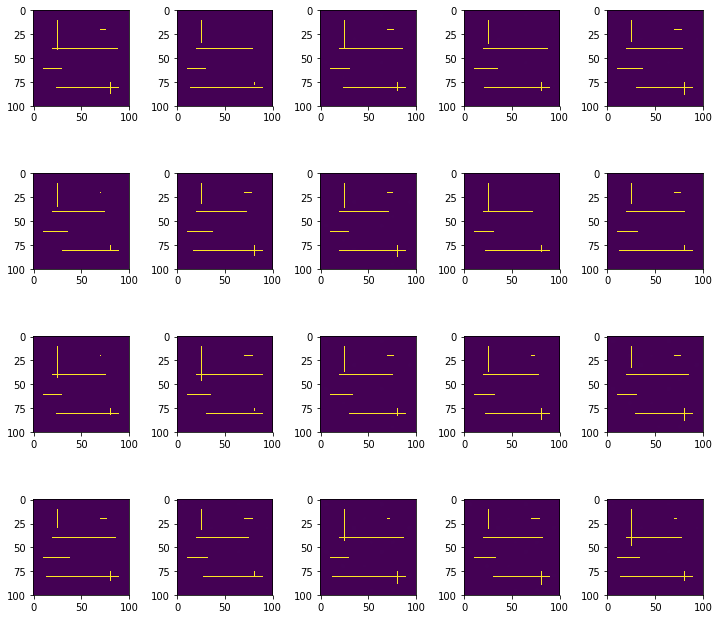}
\caption{Heterogeneous fields, the yellow strips are the channels}
\label{hfield}
\end{figure}

We consider solving (\ref{eqn1})-(\ref{eqn2}) for a heterogeneous field with moving channels and changing background. Let us denote the heterogeneous field as $\kappa(x)$, where $x\in [0,1]^2$, then $\kappa(x) = 1000$ if $x$ is in some channels which will be illustrated later and otherwise,
\begin{align*}
    \kappa(x) = e^{\eta\cdot sin(7\pi x)sin(8\pi y)+sin(10\pi x)sin(12\pi y)},
\end{align*}
where $\eta$ follows discrete uniform distribution in $[0,1]$. The channels are moving and we include cases of the intersection of two channels and formation and disappearance of the channels in the fields. In Figure~\ref{hfield}, we demonstrate $20$ heterogeneous fields. 
It can be observed from the images that, vertical channel (at around $x = 30$) (not always) intersects with horizontal channels (at around $y = 40$); and the channel at $x = 75, y = 25$ demonstrates the case of generation and degeneration of a channel. 

We train the network using 600 samples using the Adam gradient descent. We find that the cluster assignment of $600$ realizations in uncertain space is stable(fixed) when the gradient descent epoch reaches a certain number, so we set the stopping criteria to be: the assignment does not change for $100$ iteration epochs; and the maximum number of iteration epochs is set to be $1000 $. We also find that the coefficients in (\ref{loss_all}) can affect the training result. We set $\lambda_1 = \lambda_2 = \lambda_3 = 1$.

It should be noted that we train the network locally in each coarse neighborhood.
The fine mesh element has size $1/100$ and $5$ fine elements are merged into one coarse element.

\subsection{Results}
We will present the numerical results of the proposed method in this section. We are going to show the cluster assignment experiment first, followed by two other experiments which demonstrate the error of the method.
\subsubsection{Cluster Assignment in a Local Coarse Element}
\begin{figure}[H]
\centering
\includegraphics[scale = 0.4]{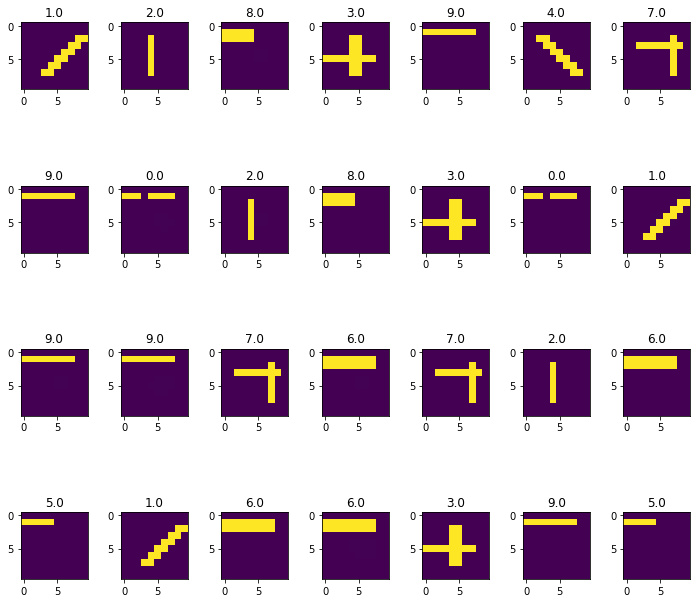}
\caption{Cluster results of 28 samples, images shown are heterogeneous fields, the number on top of each image is the cluster assignment ID number.}
\label{plot}
\end{figure}

\begin{figure}[H]
\centering
\includegraphics[scale = 0.3]{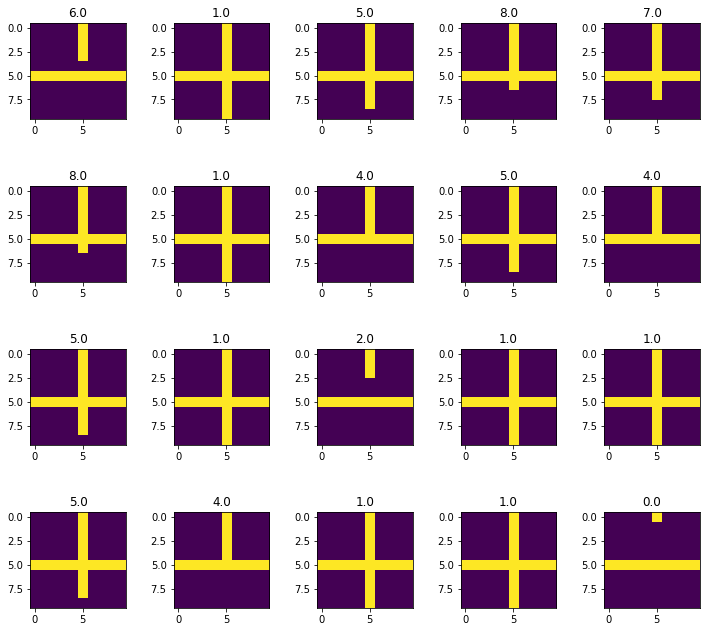}
\caption{Cluster results of 20 samples, images shown are heterogeneous fields, the number on top of each image is the cluster assignment ID number.}
\label{plot2}
\end{figure}

Before diving into the error analysis, we will show some of the cluster results in a local neighborhood. In this neighborhood, we manually created the cases such as: the extraction of a channel (longer) , the expansion of a channel(wider), the discontinuity of a channel, the diagonal channels, the intersection of channels an so on. In Figure~\ref{plot}, the number on top of each image is the cluster assignment ID number.
We also demonstrate the clustering result in Figure~\ref{plot2} of another neighborhood which is around $(25,45)$ in Figure~\ref{hfield}.
From the results in both Figure~\ref{plot} and Figure~\ref{plot2}, we observe that our proposed clustering algorithm based on deep learning
is able create a good clustering of the parameter space. That is, heterogeneous coefficients with similar spatial structures are grouped in the same cluster. 

\subsubsection{Relation of Error and the Number of Clusters}
In this section, we will demonstrate the error change when the hyperparamter --- the number of clusters increases.
Given a new realization $\kappa(x, \hat{s})$ where $\hat{s}$ denotes the parameter and a fixed neighborhood, suppose the neighborhood of this realization will be fitted into cluster $S_i$ by the model trained. We compute $\Bar{\kappa}_i = \frac{1}{|S_i|}\sum_{j = 1}^{|S_i|}\kappa_{ij}$ where $|S_i|$ is the number of points in this cluster $S_i$. 
The GMsFEM basis of this neighborhood can then be derived using $\Bar{\kappa}_i$.
We finally construct the solution using the GMsFEM basis pre-computed in all neighborhoods.
We define the $l_2$ relative error as :
\begin{align}
    ratio = \frac{\int_{D} (u-u_H)^2 dx}{\int_D u^2 dx}, \label{ratio}
\end{align}
where $u$ is the exact solution computed by finite element method with fine enough mesh and $u_H$ is the solution of the proposed method. We test the network on newly generated $300$ samples and take the average of the errors.

\begin{figure}[H]
\centering
\includegraphics[scale = 0.5]{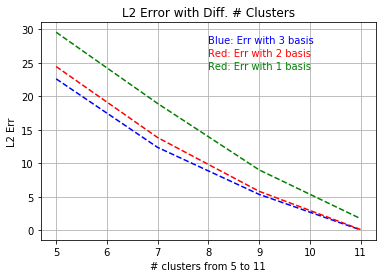}
\caption{$L_2$ error with the number of clusters changes, colors represent the number of GMsFEM basis}
\label{nbclusters}
\end{figure}

In this experiment, we calculate the $l_2$ relative error
with the number of clusters increases. The number of clusters ranges from 
$5$ to $11$; and for each case, we will train the model and compute the $l_2$ relative error.
The result can be seen in Figure~\ref{nbclusters} and it can be observed from the picture that, the error is decreasing with the number of cluster increases.

\subsubsection{Comparison of Cluster-based Method with Tradition Method}
\begin{figure}[H]
\centering
\includegraphics[scale = 0.5]{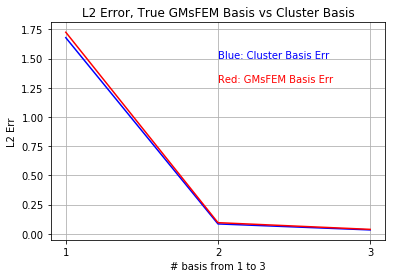}
\caption{$L_2$ error cluster solution ($11$ clusters) vs solution by real $\kappa(x, \hat{s})$. Color represents number of basis}
\label{comp}
\end{figure}

In the second experiments, we first compute the $l_2$ relative error (defined in (\ref{ratio}) with $u_H$ denotes the GMsFEM solution) of traditional GMsFEM method with given $\kappa(x, \hat{s})$. 
This means that the construct multiscale basis functions using the particular realization $\kappa(x, \hat{s})$.
We then compare this error with the cluster method proposed (11 clusters). The comparison can be seen in Figure~\ref{comp}.
It can be seen that the difference is negligible when the number of clusters reaches $11$. We can then benefit from the deep learning; i.e., the fitting of $\kappa(x, \hat{s})$ into a cluster is fast; and since we will use the pre-computed basis, we also save time on computing the GMsFEM basis.

\subsection{Effect of the Adversary Net}
The target of this task is not the learning of multiscale basis; the multiscale basis in this work is just a supervision of learning the cluster. However to demonstrate the effectiveness of the adversary network, we also test the the effect of the adversary net.
There are many hyper-parameters like the number of clusters and coefficients of the loss function which can affect the result; so to reduce the influence from the clustering, we remove the clustering loss from the training, so this is a generative task which will generate the multiscale basis from the output of the first network in Figure \ref{entirenet}.
The loss function now can be defined as:
\begin{align}
   \min_{\theta_G,\theta_F} \lambda_1 R+\lambda_2 A, \label{msbasis_ln}
\end{align}
where $R$ and $A$ are defined in (\ref{recon_err}) and (\ref{adv_err}) separately; and $\lambda_1$ and $\lambda_1$ are both set to be $1$.
We compute the relative error (\ref{ratio}) first by using the learnt multiscale basis which is trained by (\ref{msbasis_ln}); and second by using the multiscale basis trained without the adversary loss (\ref{adv_err}), i.e.,
\begin{align}
   \min_{\theta_G,\theta_F} A \label{msbasis_wo}.
\end{align}
The $l_2$ relative error improves from $41.120439$ to $36.760918$ if we add one middle layer from the adversary net. 

We also calculate the MSE difference of two learnt basis (by loss (\ref{msbasis_ln}) and (\ref{msbasis_wo}) separately) and real multiscale basis, i.e., we calculate $\| B_\text{learnt basis}-B_\text{real basis}\|_{MSE}$, where $B_\text{learnt basis}$ refers to two basis trained with (\ref{msbasis_ln}) and (\ref{msbasis_wo}) separately and $ B_\text{real basis} $ is the real multiscale basis formed using the input heterogeneous field.
The MSE amazingly decreases from $0.9073400$ to $0.748312$ if we use basis trained with the adversary loss (\ref{msbasis_ln}). This can show the benefit from the adversary net. 

\section{Conclusion}

We propose a deep learning clustering technique within GMsFEM to solve flows in heterogeneous media. The main idea is to cluster the uncertainty space such that we can reduce the number of multiscale basis functions for each coarse block across the uncertainty space.  We propose the adversary loss motivated by the perceptual loss in the computer vision task. 
We use convolutional neural networks combined with some techniques in
adversary neural networks, where the loss function is composed
of several parts that includes terms related to clusters and reconstruction
of basis functions. We present numerical results for channelized permeability
fields in the examples of flows in porous media.
In future, we would like to study the relation between 
convolutional layers and quantities related to multiscale basis functions.


\section*{Acknowledgement}

Eric Chung's work is partially supported by the Hong Kong RGC General Research Fund (Project numbers 14304719 and 14302018) and the CUHK Faculty of Science Direct Grant 2018-19.  YE would like to thank the partial support from NSF 1620318 and NSF Tripod 1934904. YE would also like to acknowledge the support of Mega-grant of the Russian Federation Government (N 14.Y26.31.0013).

\bibliographystyle{plain}

\bibliography{references_v1}

\begin{thebibliography}{10}

\bibitem{segnet}
Vijay Badrinarayanan, Alex Kendall, and Roberto Cipolla.
\newblock Segnet: A deep convolutional encoder-decoder architecture for image
  segmentation.
\newblock {\em IEEE transactions on pattern analysis and machine intelligence},
  39(12):2481--2495, 2017.

\bibitem{bishop}
Christopher~M Bishop.
\newblock {\em Pattern recognition and machine learning}.
\newblock springer, 2006.

\bibitem{gmsfem2}
Victor~M Calo, Yalchin Efendiev, Juan Galvis, and Mehdi Ghommem.
\newblock Multiscale empirical interpolation for solving nonlinear pdes.
\newblock {\em Journal of Computational Physics}, 278:204--220, 2014.

\bibitem{snapshot}
Victor~M Calo, Yalchin Efendiev, Juan Galvis, and Guanglian Li.
\newblock Randomized oversampling for generalized multiscale finite element
  methods.
\newblock {\em Multiscale Modeling \& Simulation}, 14(1):482--501, 2016.

\bibitem{cluster1}
Mathilde Caron, Piotr Bojanowski, Armand Joulin, and Matthijs Douze.
\newblock Deep clustering for unsupervised learning of visual features.
\newblock In {\em Proceedings of the European Conference on Computer Vision
  (ECCV)}, pages 132--149, 2018.

\bibitem{cheung2019deep}
Siu~Wun Cheung, Eric~T Chung, Yalchin Efendiev, Eduardo Gildin, Yating Wang,
  and Jingyan Zhang.
\newblock Deep global model reduction learning in porous media flow simulation.
\newblock {\em Computational Geosciences}, pages 1--14, 2019.

\bibitem{ElasticGMsFEM}
E.~Chung, Y.~Efendiev, and S.~Fu.
\newblock Generalized multiscale finite element method for elasticity
  equations.
\newblock {\em International Journal on Geomathematics}, 5(2):225--254, 2014.

\bibitem{ye1}
E.~T. Chung, Y.~Efendiev, and G.~Li.
\newblock An adaptive {GM}s{FEM} for high contrast flow problems.
\newblock {\em J. Comput. Phys.}, 273:54--76, 2014.

\bibitem{gmsfem9}
Eric Chung, Yalchin Efendiev, and Thomas~Y Hou.
\newblock Adaptive multiscale model reduction with generalized multiscale
  finite element methods.
\newblock {\em Journal of Computational Physics}, 320:69--95, 2016.

\bibitem{chung2017DGstokes}
Eric Chung, Maria Vasilyeva, and Yating Wang.
\newblock A conservative local multiscale model reduction technique for stokes
  flows in heterogeneous perforated domains.
\newblock {\em Journal of Computational and Applied Mathematics}, 321:389--405,
  2017.

\bibitem{gmsfem3}
Eric~T Chung, Yalchin Efendiev, and Wing~Tat Leung.
\newblock Residual-driven online generalized multiscale finite element methods.
\newblock {\em Journal of Computational Physics}, 302:176--190, 2015.

\bibitem{gmsfem4}
Eric~T Chung, Yalchin Efendiev, and Wing~Tat Leung.
\newblock An online generalized multiscale discontinuous galerkin method
  (gmsdgm) for flows in heterogeneous media.
\newblock {\em Communications in Computational Physics}, 21(2):401--422, 2017.

\bibitem{paper}
Eric~T Chung, Yalchin Efendiev, Wing~Tat Leung, and Zhiwen Zhang.
\newblock Cluster-based generalized multiscale finite element method for
  elliptic pdes with random coefficients.
\newblock {\em Journal of Computational Physics}, 371:606--617, 2018.

\bibitem{gmsfem5}
Eric~T Chung, Yalchin Efendiev, and Guanglian Li.
\newblock An adaptive gmsfem for high-contrast flow problems.
\newblock {\em Journal of Computational Physics}, 273:54--76, 2014.

\bibitem{gmsfem6}
Eric~T Chung, Yalchin Efendiev, Guanglian Li, and Maria Vasilyeva.
\newblock Generalized multiscale finite element methods for problems in
  perforated heterogeneous domains.
\newblock {\em Applicable Analysis}, 95(10):2254--2279, 2016.

\bibitem{dong}
Chao Dong, Chen~Change Loy, Kaiming He, and Xiaoou Tang.
\newblock Image super-resolution using deep convolutional networks.
\newblock {\em IEEE transactions on pattern analysis and machine intelligence},
  38(2):295--307, 2015.

\bibitem{gmsfem1}
Yalchin Efendiev, Juan Galvis, and Thomas~Y Hou.
\newblock Generalized multiscale finite element methods (gmsfem).
\newblock {\em Journal of Computational Physics}, 251:116--135, 2013.

\bibitem{gmsfem7}
Yalchin Efendiev, Juan Galvis, R~Lazarov, M~Moon, and Marcus Sarkis.
\newblock Generalized multiscale finite element method. symmetric interior
  penalty coupling.
\newblock {\em Journal of Computational Physics}, 255:1--15, 2013.

\bibitem{gmsfem8}
Yalchin Efendiev, Juan Galvis, Guanglian Li, and Michael Presho.
\newblock Generalized multiscale finite element methods: Oversampling
  strategies.
\newblock {\em International Journal for Multiscale Computational Engineering},
  12(6), 2014.

\bibitem{goodfellow}
Ian Goodfellow, Yoshua Bengio, and Aaron Courville.
\newblock {\em Deep learning}.
\newblock MIT press, 2016.

\bibitem{gan}
Ian Goodfellow, Jean Pouget-Abadie, Mehdi Mirza, Bing Xu, David Warde-Farley,
  Sherjil Ozair, Aaron Courville, and Yoshua Bengio.
\newblock Generative adversarial nets.
\newblock In {\em Advances in neural information processing systems}, pages
  2672--2680, 2014.

\bibitem{msfem1}
Thomas~Y Hou and Xiao-Hui Wu.
\newblock A multiscale finite element method for elliptic problems in composite
  materials and porous media.
\newblock {\em Journal of computational physics}, 134(1):169--189, 1997.

\bibitem{senet}
Jie Hu, Li~Shen, and Gang Sun.
\newblock Squeeze-and-excitation networks.
\newblock In {\em Proceedings of the IEEE conference on computer vision and
  pattern recognition}, pages 7132--7141, 2018.

\bibitem{dense}
Gao Huang, Zhuang Liu, Laurens Van Der~Maaten, and Kilian~Q Weinberger.
\newblock Densely connected convolutional networks.
\newblock In {\em Proceedings of the IEEE conference on computer vision and
  pattern recognition}, pages 4700--4708, 2017.

\bibitem{condition}
Phillip Isola, Jun-Yan Zhu, Tinghui Zhou, and Alexei~A Efros.
\newblock Image-to-image translation with conditional adversarial networks.
\newblock In {\em Proceedings of the IEEE conference on computer vision and
  pattern recognition}, pages 1125--1134, 2017.

\bibitem{msfem2}
Patrick Jenny, SH~Lee, and Hamdi~A Tchelepi.
\newblock Multi-scale finite-volume method for elliptic problems in subsurface
  flow simulation.
\newblock {\em Journal of Computational Physics}, 187(1):47--67, 2003.

\bibitem{perceptual}
Justin Johnson, Alexandre Alahi, and Li~Fei-Fei.
\newblock Perceptual losses for real-time style transfer and super-resolution.
\newblock In {\em European conference on computer vision}, pages 694--711.
  Springer, 2016.

\bibitem{kl1}
Kari Karhunen.
\newblock {\em {\"U}ber lineare Methoden in der Wahrscheinlichkeitsrechnung},
  volume~37.
\newblock Sana, 1947.

\bibitem{kim1}
Jiwon Kim, Jung Kwon~Lee, and Kyoung Mu~Lee.
\newblock Accurate image super-resolution using very deep convolutional
  networks.
\newblock In {\em Proceedings of the IEEE conference on computer vision and
  pattern recognition}, pages 1646--1654, 2016.

\bibitem{laplacian}
Wei-Sheng Lai, Jia-Bin Huang, Narendra Ahuja, and Ming-Hsuan Yang.
\newblock Deep laplacian pyramid networks for fast and accurate
  super-resolution.
\newblock In {\em Proceedings of the IEEE conference on computer vision and
  pattern recognition}, pages 624--632, 2017.

\bibitem{supergan}
Christian Ledig, Lucas Theis, Ferenc Husz{\'a}r, Jose Caballero, Andrew
  Cunningham, Alejandro Acosta, Andrew Aitken, Alykhan Tejani, Johannes Totz,
  Zehan Wang, et~al.
\newblock Photo-realistic single image super-resolution using a generative
  adversarial network.
\newblock In {\em Proceedings of the IEEE conference on computer vision and
  pattern recognition}, pages 4681--4690, 2017.

\bibitem{enhanced}
Bee Lim, Sanghyun Son, Heewon Kim, Seungjun Nah, and Kyoung Mu~Lee.
\newblock Enhanced deep residual networks for single image super-resolution.
\newblock In {\em Proceedings of the IEEE conference on computer vision and
  pattern recognition workshops}, pages 136--144, 2017.

\bibitem{kmeans}
Stuart Lloyd.
\newblock Least squares quantization in pcm.
\newblock {\em IEEE transactions on information theory}, 28(2):129--137, 1982.

\bibitem{long}
Jonathan Long, Evan Shelhamer, and Trevor Darrell.
\newblock Fully convolutional networks for semantic segmentation.
\newblock In {\em Proceedings of the IEEE conference on computer vision and
  pattern recognition}, pages 3431--3440, 2015.

\bibitem{deconvolution}
Hyeonwoo Noh, Seunghoon Hong, and Bohyung Han.
\newblock Learning deconvolution network for semantic segmentation.
\newblock In {\em Proceedings of the IEEE international conference on computer
  vision}, pages 1520--1528, 2015.

\bibitem{simonyan2014very}
Karen Simonyan and Andrew Zisserman.
\newblock Very deep convolutional networks for large-scale image recognition.
\newblock {\em arXiv preprint arXiv:1409.1556}, 2014.

\bibitem{inception}
Christian Szegedy, Sergey Ioffe, Vincent Vanhoucke, and Alexander~A Alemi.
\newblock Inception-v4, inception-resnet and the impact of residual connections
  on learning.
\newblock In {\em Thirty-first AAAI conference on artificial intelligence},
  2017.

\bibitem{rethinking}
Christian Szegedy, Vincent Vanhoucke, Sergey Ioffe, Jon Shlens, and Zbigniew
  Wojna.
\newblock Rethinking the inception architecture for computer vision.
\newblock In {\em Proceedings of the IEEE conference on computer vision and
  pattern recognition}, pages 2818--2826, 2016.

\bibitem{recursive}
Ying Tai, Jian Yang, and Xiaoming Liu.
\newblock Image super-resolution via deep recursive residual network.
\newblock In {\em Proceedings of the IEEE conference on computer vision and
  pattern recognition}, pages 3147--3155, 2017.

\bibitem{memnet}
Ying Tai, Jian Yang, Xiaoming Liu, and Chunyan Xu.
\newblock Memnet: A persistent memory network for image restoration.
\newblock In {\em Proceedings of the IEEE international conference on computer
  vision}, pages 4539--4547, 2017.

\bibitem{vasilyeva2019learning}
Maria Vasilyeva, Wing~T Leung, Eric~T Chung, Yalchin Efendiev, and Mary
  Wheeler.
\newblock Learning macroscopic parameters in nonlinear multiscale simulations
  using nonlocal multicontinua upscaling techniques.
\newblock {\em JCP, accepted, arXiv preprint arXiv:1907.02921}, 2019.

\bibitem{wang2019prediction}
Min Wang, Siu~Wun Cheung, Eric~T Chung, Yalchin Efendiev, Wing~Tat Leung, and
  Yating Wang.
\newblock Prediction of discretization of gmsfem using deep learning.
\newblock {\em Mathematics}, 7(5):412, 2019.

\bibitem{wang2020reduced}
Min Wang, Siu~Wun Cheung, Wing~Tat Leung, Eric~T Chung, Yalchin Efendiev, and
  Mary Wheeler.
\newblock Reduced-order deep learning for flow dynamics. the interplay between
  deep learning and model reduction.
\newblock {\em Journal of Computational Physics}, 401:108939, 2020.

\bibitem{wang2020deep}
Yating Wang, Siu~Wun Cheung, Eric~T Chung, Yalchin Efendiev, and Min Wang.
\newblock Deep multiscale model learning.
\newblock {\em Journal of Computational Physics}, 406:109071, 2020.

\bibitem{cluster2}
Junyuan Xie, Ross Girshick, and Ali Farhadi.
\newblock Unsupervised deep embedding for clustering analysis.
\newblock In {\em International conference on machine learning}, pages
  478--487, 2016.

\bibitem{cluster3}
Bo~Yang, Xiao Fu, Nicholas~D Sidiropoulos, and Mingyi Hong.
\newblock Towards k-means-friendly spaces: Simultaneous deep learning and
  clustering.
\newblock In {\em Proceedings of the 34th International Conference on Machine
  Learning-Volume 70}, pages 3861--3870. JMLR. org, 2017.

\bibitem{selfattentiongan}
Han Zhang, Ian Goodfellow, Dimitris Metaxas, and Augustus Odena.
\newblock Self-attention generative adversarial networks.
\newblock {\em arXiv preprint arXiv:1805.08318}, 2018.

\bibitem{residualchannel}
Yulun Zhang, Kunpeng Li, Kai Li, Lichen Wang, Bineng Zhong, and Yun Fu.
\newblock Image super-resolution using very deep residual channel attention
  networks.
\newblock In {\em Proceedings of the European Conference on Computer Vision
  (ECCV)}, pages 286--301, 2018.

\bibitem{residualdense}
Yulun Zhang, Yapeng Tian, Yu~Kong, Bineng Zhong, and Yun Fu.
\newblock Residual dense network for image super-resolution.
\newblock In {\em Proceedings of the IEEE conference on computer vision and
  pattern recognition}, pages 2472--2481, 2018.

\end{thebibliography}

\end{document}